\documentclass[12pt]{amsart}

\newcommand{\B}{\mathcal{B}}
\newcommand{\BG}{\B_\Gamma}
\newcommand{\BL}{\B_\Lambda}
\newcommand{\BT}{\B_\Tau}

\newcommand{\BO}{\B_\Omega}
\newcommand{\CG}{C_\Gamma}
\newcommand{\CL}{C_\Lambda}
\newcommand{\CT}{C_\Tau}

\newcommand{\CO}{C_\Omega}
\newcommand{\sone}{\mathsf{S}_1}    \newcommand{\sfin}{\mathsf{S}_{fin}}
\newcommand{\ufin}{\mathsf{U}_{fin}}
    
\newcommand{\seq}[1]{\{#1\}_{n\in\N}}
\newcommand{\Union}{\bigcup}
\newcommand{\nin}{\not\in}
\newcommand{\cF}{\mathcal{F}}
\newcommand{\cU}{\mathcal{U}}
\newcommand{\cV}{\mathcal{V}}

\newcommand{\fU}{\mathfrak{U}}
\newcommand{\fV}{\mathfrak{V}}

\newcommand{\naturals}{{\mathbb N}}
\newcommand{\N}{\naturals}

\newcommand{\sbst}{\subseteq}
\newcommand{\by}[2]{\hfill\emph{#1}, #2}
\newcommand{\Tau}{\mathrm{T}}

\newcommand{\be}{\begin{enumerate}}
\newcommand{\ee}{\end{enumerate}}
\newcommand{\bi}{\begin{itemize}}
\newcommand{\ei}{\end{itemize}}
\renewcommand{\i}{\item}
\newcommand{\SPMBul}{\textbf{$\mathcal{SPM}$ BULLETIN}}
\newcommand{\BulEnd}{\par\bigskip\noindent
Boaz Tsaban\\
\emph{E-mail}: tsaban@math.huji.ac.il\\
\emph{URL}: http://www.cs.biu.ac.il/\~{}tsaban}

\title[\SPMBul{} \textbf{1} (Jan 2003)]{%
\SPMBul\\[0.5cm]
Issue number 1: January 2003 \textsc{ce}}

\begin{document}
\maketitle

\tableofcontents

\section{Introduction}

We are glad to introduce the first semi-formal issue of the
\SPMBul{}. This bulletin began as e-mails sent to a small list of
colleagues working, or interested, in the field of
\emph{Selection Principles in Mathematics} (SPM).
The original mailing list consisted of attendants in the
June 2002 Lecce workshop on \emph{Coverings, Selesctions, and Games in Topology},
and then we decided to try to extend it to include
all mathematicians which made or make a significant contribution to the field,
as well as some good mathematicians which are interested in the field.
Currently, we are very far from this utopian goal, as the current list contains
only about $40$ recipients (these are the ones which the editor knows).
We hope that the recipients will be able to supply more names
(and e-mails, to the address at the end of this issue; note that
we do not include new members in our list before we ask them if they
would like to join us).

\subsection{What is SPM?}\label{whatisSPM}\footnote{The
following is arguable and, as everything else in the \SPMBul{}, reflects
only the personal opinion of the editor.}
For a good answer to this we are afraid that the reader should wait until
Marion Scheepers will finish the writing of
his introductory paper on the field and make it publicly available.
Meanwhile, we will say vaguely that this field is about all sorts of studies
of diagonalization arguments, especially in topology (covering properties,
sequences of covers, etc.) and infinite combinatorics (cardinal characteristics
of the continuum), and their applications to other areas of mathematics
(function spaces, game theory, group theory, etc.).
Some readers know this field in other names, e.g.:
\bi
\i \emph{Topological Games} -- Telgarski, et.\ al.,
\i \emph{Combinatorics of open/Borel covers} -- Scheepers,
\i \emph{Topological Diagonalizations} -- folklore; and
\i \emph{Infinite-combinatorial topology} -- Tsaban (?).
\ei
Some other are better acquainted with the more classical occurrences of
this field, which are of more particular-case nature, e.g.,
Menger and Hurewicz properties, Rothberger property $C''$,
Gerlits-Nagy $\gamma$-sets, etc.
Certainly, we view the field of
\emph{Special (or: Totally Imperfect) Sets of Reals} as a subfield of
this field, or at least as one with large overlap with SPM.

\subsection{Aims}
\subsubsection{New results}
The \SPMBul{} is intended to be a quick and informal mean to transmit knowledge
between mathematicians working in the field. Being purely electronic, and having its
focused audience, we hope to be among the first to report new results in the field.
We will also announce here new papers appearing in the \emph{mathematics ArXiv}
which are relevant to SPM, and quote their abstracts.

\subsubsection{Discussions}
Another goal is to discuss various open problems in the field.
Problems presented here do not have to be very difficult -- it may be the case that
the one presenting them overlooked something or just does not know of a result which
implies an answer. The solutions to problems (if found) and related suggestions
will be announced in following issues.
We also hope to have a section called \emph{Problem of the month} which will include
a really tough and important problem.

\subsection{Contributions}
It is expected that most of the contributions (announcements, discussions, and open problems)
will come from  the recipients of this
bulletin. The contributor's name will appear at the end of his or her contribution.
Contributions should be made in \LaTeX{} (or otherwise in some other \TeX{} format
or plain text), and emailed to the editor.
The authors are urged to use as standard notation as possible, or otherwise give
a reference to where the notation is explained.

Clearly contributions to this bulletin would not require any transfer of copyright,
and material presented informally here can be edited and published anywhere the
author likes.

\subsection{Publication and citation}
The \SPMBul{} will appear \emph{at most} monthly, so not to overload the
e-mailboxes of the recipients.
It will also be published electronically on
the \emph{arXiv.org e-Print archive}
\begin{quote}
http://arxiv.org/
\end{quote}
(in the categories: General Topology, Logic, and Combinatorics).

The submission to the \emph{ArXiv} may occur later than the reception
via e-mail. However, it is useful for citations in papers, e.g., when
solving a problem posed in the \SPMBul{} one can give a reference in the
form:\\[0.1cm]
Giorgy Giorge, \emph{Does $\binom{\Omega}{\Tau}=\binom{\Omega}{\Gamma}$?},
SPM Bulletin \textbf{1} (2003), 2--3,\\
\verb|http://arxiv.org/abs/math/0123456|

\subsection{Help needed}
In addition to mathematical contributions, any suggestions with regard to the form and
content of this bulletin will be appreciated.
Also, it will be impossible not to miscredit anyone when announcing results, especially due
to the rapid character of this bulletin.
Please send us the corrections and we will publish them in subsequent issues.
(We will also appreciate any grammar and other corrections, but these would not be published
unless this is really important.)
Finally, we are totally ignorant about legal issues, and would appreciate any suggestion concerning
these.

\subsection{Subscription and unsubscription}
Anyone who wishes to receive this bulletin (free of charge!)
to his e-mailbox should contact the editor
at the address given at the end of this issue.
Unsubscription is also done via e-mailing the editor.

\section{Some basic notation}\label{basedefs}
We introduce here the most basic standard notations in the field.
\subsection{Selection principles}
Let $\fU$ and $\fV$ be collections of covers of a space $X$.
Following are selection hypotheses which $X$ might satisfy or not
satisfy.\footnote{Often these hypotheses are identified with the class of
all spaces satisfying them.}
\begin{itemize}
\item[$\sone(\fU,\fV)$:]
For each sequence $\seq{\cU_n}$ of members of $\fU$,
there is a sequence
$\seq{V_n}$ such that for each $n$ $V_n\in\cU_n$, and $\seq{V_n}\in\fV$.
\item[$\sfin(\fU,\fV)$:]
For each sequence $\seq{\cU_n}$
of members of $\fU$, there is a sequence
$\seq{\cF_n}$ such that each $\cF_n$ is a finite
(possibly empty) subset of $\cU_n$, and
$\Union_{n\in\N}\cF_n\in\fV$.
\item[$\ufin(\fU,\fV)$:]
For each sequence
$\seq{\cU_n}$ of members of $\fU$
\emph{which do not contain a finite subcover},
there exists a sequence $\seq{\cF_n}$
such that for each $n$ $\cF_n$ is a finite (possibly empty) subset of
$\cU_n$, and
$\seq{\cup\cF_n}\in\fV$.
\end{itemize}
Also define the property ``\emph{$\fU$ choose $\fV$}'' as follows.
\bi
\i[$\binom{\fU}{\fV}$:]
For each $\cU\in\fU$ there exists $\cV\sbst\cU$ such that $\cV\in\fV$.
\ei
The hypotheses $\sone(\fU,\fV)$, $\sfin(\fU,\fV)$ and $\binom{\fU}{\fV}$
make sense in any context where $\fU$ and $\fV$ are
collections of subsets of an infinite set.

\subsection{Thick covers}
Let $\cU$ be a collection of subsets of $X$ such that $X$ is not
contained in any member of $\cU$. $\cU$ is:
\be
\i A \emph{large cover} of $X$ if each $x\in X$ is contained in infinitely many members of $\cU$,
\i An \emph{$\omega$-cover} of $X$ if each finite subset of $X$ is contained in some member of $\cU$,
\i A \emph{$\tau$-cover} of $X$ if it is a large cover of $X$,
and for each $x,y\in X$, either $\{U\in\cU : x\in U, y\nin U\}$ is finite, or
$\{U\in\cU : y\in U, x\nin U\}$ is finite; and
\i A $\gamma$-cover of $X$ if $\cU$ is infinite, and
each $x\in X$ belongs to all but finitely many
members of $\cU$.
\ee
Let $\Lambda$, $\Omega$, $\Tau$, and $\Gamma$
denote the collections of open large covers, $\omega$-covers,
$\tau$-covers, and
$\gamma$-covers of $X$, respectively. Also,
let $\BL,\BO,\BT,\BG$ (respectively, $\CL,\CO,\CT,\CG$)
be the corresponding \emph{countable Borel} (respectively, \emph{clopen}) covers of $X$.

\section{Announcements}

\subsection{A Mad $Q$-set}
A MAD (maximal almost disjoint) family is an infinite subset $A$ of the
infinite subsets of $\{0,1,2,\dots\}$ such that any two elements of $A$ intersect in a
finite set and every infinite subset of $\{0,1,2,\dots\}$ meets some element of $A$
in an infinite set. A $Q$-set is an uncountable set of reals such that every
subset is a relative $G_\delta$ set. It is shown that it is relatively consistent
with ZFC that there exists a MAD family which is also a $Q$-set in the topology
in inherits a subset of the Power set of $\{0,1,2,\dots\}$, i.e.\ the Cantor set.

\verb|http://arXiv.org/abs/math/0212335|

\by{Arnold W.\ Miller}{miller@math.wisc.edu}

\subsection{On $\lambda'$-sets}
A subset $X$ of the Cantor space, $2^\omega$, is a $\lambda'$-set iff for every
countable subset $Y$ of the Cantor space $Y$ is relatively $G_\delta$ in $X$ union $Y$. In
this paper we prove two forcing results about $\lambda'$-sets. First we show
that it is consistent that every $\lambda'$-set is a $\gamma$-set. Secondly we
show that is independent whether or not every $\dagger-\lambda$-set is a
$\lambda'$-set.

\verb|http://arXiv.org/abs/math/0212336|

\by{Arnold W.\ Miller}{miller@math.wisc.edu}

\subsection{Problem paper}
In the Lecce Workshop on \emph{Coverings, Selesctions, and Games in Topology}
(June 2002), it was decided that we write a paper on open problems in the field
(see Section \ref{whatisSPM} above).
It is hoped that this paper will be a useful source of problems for professionals
as well as beginners in the field, and this is a wonderful way to make the problems
you are interested in widely available.

Please submit to me a short description of your problems together with definitions
(if not standard in SPM), motivation, and possible tools or references which may
be helpful in dealing with these problems. If the problem was stated earlier in
some public form, a reference can be helpful.

Please submit only the problems which you consider more important among the problems
you know of.

\by{Boaz Tsaban}{tsaban@math.huji.ac.il}

\section{Problem of the month}
Using the notation given in \ref{basedefs}, the problem can be stated easily.
$\binom{\Omega}{\Gamma}$ is the celebrated Gerlits-Nagy $\gamma$-property \cite{GN}.
As $\Gamma\sbst\Tau$, this property implies $\binom{\Omega}{\Tau}$.
Thus far, all examples of sets not satisfying $\binom{\Omega}{\Gamma}$
turned out not to satisfy $\binom{\Omega}{\Tau}$.
\begin{quote}
Are the properties $\binom{\Omega}{\Gamma}$ and $\binom{\Omega}{\Tau}$
(provably) equivalent for zero-dimensional sets of reals $X$?
\end{quote}
This problem seem to rise almost in every study of $\tau$-covers, e.g.,
\cite{tautau, strongdiags, splittability}, and a positive answer would solve several
open problems in the field. A negative answer should also imply (through a
bit finer analysis) a solution to several open problems.
(The last three cited references are available from the author.)

\by{Boaz Tsaban}{tsaban@math.huji.ac.il}

\BulEnd

\end{document}